\documentclass{article}[12pt]

\usepackage{verbatim,graphicx,xcolor,soul,color,amsmath,amsthm,placeins,epstopdf,enumerate}

\newtheorem*{claim}{Claim}

\title{The Mathematics of String Art Nets}
\author{Chaz Lebouthillier and Mateja \v{S}ajna \\ University of Ottawa}

\renewcommand{\phi}{{\varphi}}
\renewcommand{\a}{{\rm area}}

\begin{document}
\maketitle
\baselineskip 17pt

\begin{figure}[h!]
\centering
\includegraphics[scale=0.35]{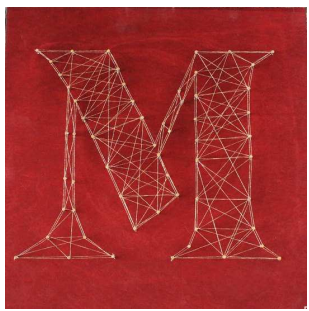}
\caption{An example of string art.}
\label{fig:M-art}
\end{figure}

\section{Introduction}

\begin{figure}[t]
\centering
\includegraphics[scale=0.35]{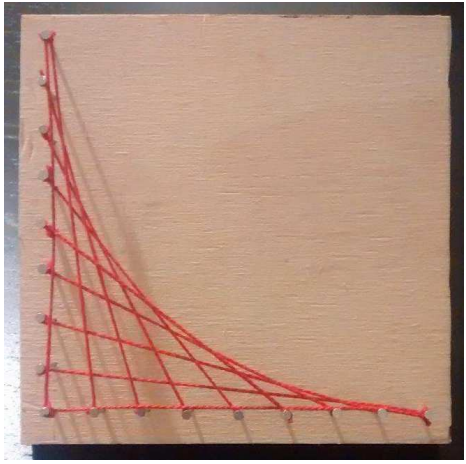}
\caption{An example of string art studied in this article.}
\label{fig:art}
\end{figure}

{\em String art} is an arrangement of pegs on a board with thread strung between these pegs to form beautiful geometric and other patterns, for example, as in  Figure~\ref{fig:M-art}. In this article, we shall consider a particular simple form of string art shown in Figure~\ref{fig:art} and its variants. Here, pegs are placed on two diverging axes, and segments of string join the first peg on one axis to the last peg on the second axis, the second peg on the first axis to the second-to-last peg on the second axis, and so on. The resulting pattern is a peculiarly shaped net. Its upper-right edge is a piecewise-linear curve whose envelope, discussed by Gregory Quenell in \cite{Que}, is a parabola. In this article, however, we would like to show that the net itself has equally fascinating mathematical properties.

\section{Prerequisites}

In this article, we shall use the following formulas. The equation of a straight line through points $(x_1,y_1)$ and $(x_2,y_2)$ is
\begin{eqnarray}
y &=& \frac{y_2 - y_1}{x_2 - x_1}\left(x - x_1\right) + y_1. \label{line}
\end{eqnarray}
The area of the triangle in the cartesian plane with vertices $P_1=(x_1,y_1)$,  $P_2=(x_2,y_2)$, and $P_3=(x_3,y_3)$ is
\begin{eqnarray}
\a ( \triangle P_1P_2P_3) &=& \frac{1}{2}\Big\vert x_1(y_2-y_3)+x_2(y_3-y_1)+x_3(y_1-y_2)\Big\vert \label{area}.
\end{eqnarray}
The sum of the first $n$ positive integers is
\begin{eqnarray}
\sum_{i=1}^n i &=& \frac{1}{2}n(n+1), \label{sum}
\end{eqnarray}
and the sum of the squares of the first $n$ positive integers is
\begin{eqnarray}
\sum_{i=1}^n i^2=\frac{1}{6}n(n+1)(2n+1). \label{squares}
\end{eqnarray}
Throughout this paper, $N$ will denote a fixed positive integer (one less than the number of pegs on one axis) and  $I=\{0,1,\ldots,N\}$.

\section{String art nets with a right-angled frame and equidistant intercepts}

\begin{figure}[t]
\centering
\includegraphics[scale=0.7]{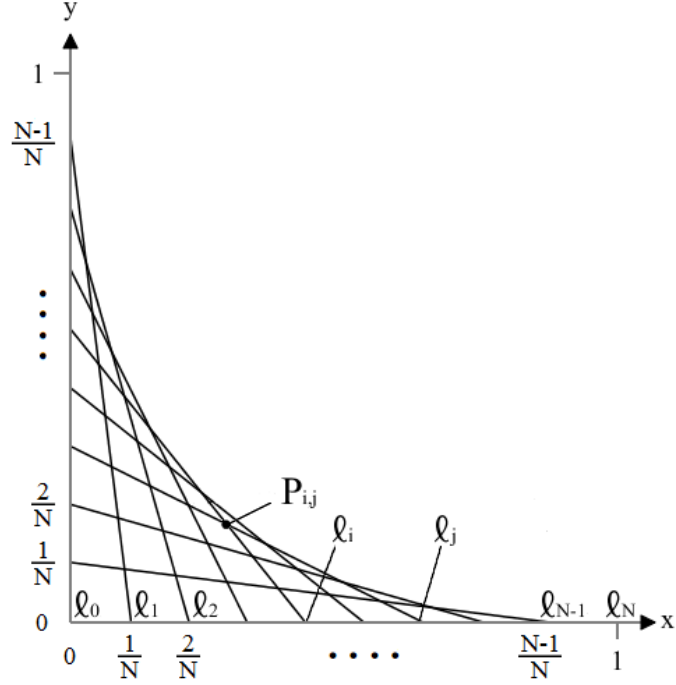}
\caption{A string art net with a right-angled frame and equidistant intercepts.}
\label{fig:net}
\end{figure}

We shall first examine the mathematical properties of a piece of string art  illustrated in Figure~\ref{fig:art}; that is, with a  right-angled frame and equidistant pegs. Our net consists of straight lines $\ell_0,\ell_1,\ldots,\ell_N$ in the first quadrant of the coordinate system (Figure~\ref{fig:net}).
For each $i \in I$, line $\ell_i$ is defined by its intercepts
\begin{eqnarray}
& X_i = \left(\frac{i}{N},0\right) \qquad \mbox{and} \qquad & Y_i = \left(0,1-\frac{i}{N}\right). \label{intercepts-basic}
\end{eqnarray}
Using Formula~(\ref{line}), the equation of line $\ell_i$ is
$$y = \frac{(1 - \frac{i}{N}) - 0}{0-\frac{i}{N}}\left(x - \frac{i}{N}\right) + 0. $$
Rearranging the terms, we obtain the canonical form
$$y = -\frac{N - i}{i}x + \Big(1 - \frac{i}{N}\Big). $$
For any $i,j \in I$, $i \ne j$, the $x$-coordinate of the intersection point of lines $\ell_i$ and $\ell_j$ is the solution to the equation
$$-\frac{N - i}{i}x + \Big(1 - \frac{i}{N}\Big) = -\frac{N - j}{j}x + \Big(1 - \frac{j}{N}\Big),$$
that is,
$$ x = \frac{1}{N^2}ij.$$
It follows that the $y$-coordinate of the intersection point is
$$y = -\frac{N - i}{i} \cdot \frac{ij}{N^2} + \Big(1 - \frac{i}{N}\Big) = \frac{1}{N^2}(N - i)(N - j).$$
Therefore, the intersection of lines  $\ell_i$ and $\ell_j$ (for $i \ne j$) is the point
\begin{eqnarray}
P_{i,j}=\frac{1}{N^2}\Big( ij,(N - i)(N - j) \Big). \label{vertex}
\end{eqnarray}

\begin{figure}[t]
\centering
\includegraphics[scale=0.6]{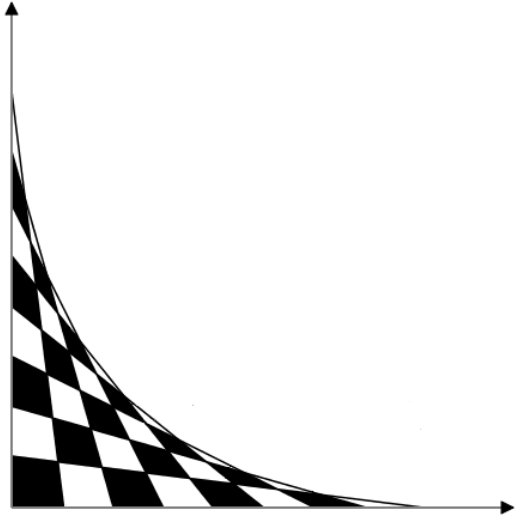}
\caption{A checkered net to better visualize the quadrilaterals.}
\label{fig:checkerednet}
\end{figure}

The lines $\ell_0,\ldots,\ell_N$ partition the area in the first quadrant (under at least one of these lines) into quadrilaterals and triangles; see Figures~\ref{fig:net} and \ref{fig:checkerednet}. As we shall see later, the areas of these quadrilaterals and triangles exhibit some surprising properties. But first, we shall determine the lengths of the sides of these shapes.

\begin{figure}[t]
\centering
\includegraphics[scale=0.45]{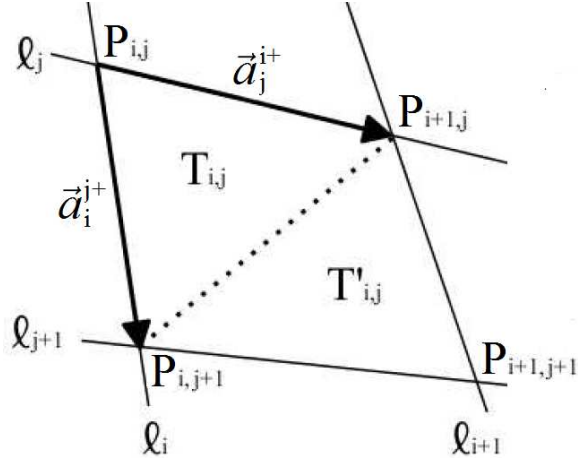}
\caption{Quadrilateral $Q_{i,j}$, and triangles $T_{i,j}$ and $T'_{i,j}$.}
\label{fig:quadri1}
\end{figure}

\begin{figure}[t]
\centering
\includegraphics[scale=0.6]{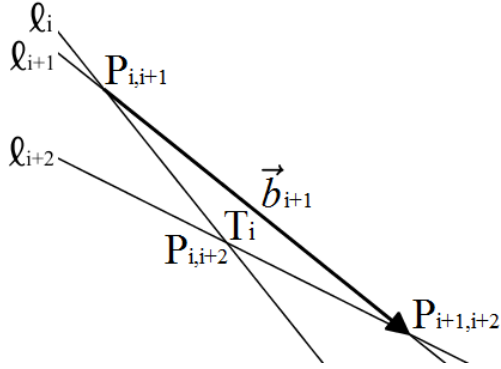}
\caption{Triangle $T_{i}$, bound by $\ell_i$,$\ell_{i+1}$, and $\ell_{i+2}$.}
\label{fig:triangle}
\end{figure}

For all $i,j \in I$ such that $j \not\in \{i-1,i,N\}$, define a vector
\begin{eqnarray}
\vec{a}_{i}^{j+} &=& P_{i,j+1} - P_{i,j} = \frac{1}{N^2}\Big( i, -(N - i) \Big), \label{a}
\end{eqnarray}
and for all $i \in \{1,2,\ldots, N-1\}$, define a vector
\begin{eqnarray}
\vec{b}_{i} = P_{i,i+1} - P_{i-1,i} = \frac{2}{N^2}\Big( i, -(N - i) \Big). \label{b}
\end{eqnarray}
Observe that $\vec{a}_{i}^{j+}$ and $\vec{b}_{i}$ are both parallel to line $\ell_i$. Vector $\vec{a}_{i}^{j+}$ has as initial point the intersection of $\ell_i$ with $\ell_j$, and as terminal point the intersection of $\ell_i$ with $\ell_{j+1}$ (Figure~\ref{fig:quadri1}), while vector $\vec{b}_{i}$ extends from the intersection of line $\ell_{i}$ with line $\ell_{i-1}$ to the intersection of line $\ell_{i}$ with line $\ell_{i+1}$ (Figure~\ref{fig:triangle}).

Fix any $i \in I$ such that $1 \le i \le N-1$, and consider the line $\ell_i$. The intersection points $P_{i,j}$, for all $j \in I$, $j \ne i$, divide the line $\ell_i$ in the first quadrant into the line segments corresponding to vectors
$$\vec{a}_{i}^{0+}, \vec{a}_{i}^{1+},\ldots, \vec{a}_{i}^{(i-2)+},\vec{b}_{i},\vec{a}_{i}^{(i+1)+},\ldots, \vec{a}_{i}^{(N-1)+}.$$
Expressions (\ref{a}) and (\ref{b}), however, show that these vectors depend only on the parameter $i$; moreover,
$$\vec{a}_{i}^{0+}= \vec{a}_{i}^{1+}=\ldots= \vec{a}_{i}^{(i-2)+} =\frac{1}{2}\vec{b}_{i}= \vec{a}_{i}^{(i+1)+}=\ldots= \vec{a}_{i}^{(N-1)+}.$$
Hence the line segments along $\ell_i$ are all of the same length, with the exception of the line segment $P_{i,i-1}P_{i,i+1}$, which is twice as long as all others!

Let's say a few more words about the length of these segments, which we shall denote by $s_i$; that is, $s_i=|\vec{a}_{i}^{j+}|$ for all $j \not\in \{ i-1,i,N\}$.
Since line $\ell_i$  is partitioned in the first quadrant into $N-2$ segments of length $s_i$ and one segment of length $2s_i$, we have
$$s_i=\frac{1}{N}\sqrt{i^2+(N-i)^2}.$$
Consequently, if $s_i=s_j$ for some $i \ne j$, then necessarily $j=N-i$. We shall refer to this observation shortly.

Next, we turn our attention to the quadrilaterals. Denote the quadrilateral with vertices $P_{i,j}, P_{i+1,j}$,  $P_{i,j+1}$, and $P_{i+1,j+1}$ by $Q_{i,j}$ (Figure~\ref{fig:quadri1}). Note that we may assume $0 \le i<i+1<j<j+1\le N$, and hence also that $N \ge 4$. We say that quadrilaterals $Q_{0,d}, Q_{1,1+d},Q_{2,2+d},\ldots,Q_{N-d-1,N-1}$, for any $d \in \{ 2,3,\ldots,N-1\}$,  lie on the same {\em diagonal} of the string art net. Indeed, in Figure~\ref{fig:checkerednet}, these quadrilaterals are of the same colour and appear to be strung on a curve that runs from the upper-left side to the bottom-right side of the net.

We first show that when a quadrilateral $Q_{i,j}$ is divided from its bottom left corner to its top right corner, the resulting triangles have the same area.
Let $T_{i,j}$ denote the triangle with vertices $P_{i,j}, P_{i+1,j}$, and $P_{i,j+1}$, and $T'_{i,j}$ the triangle with vertices $P_{i+1,j}, P_{i,j+1}$, and $P_{i+1,j+1}$ (Figure~\ref{fig:quadri1}). Using Formulas (\ref{area}) and  (\ref{vertex}), the area of triangle $T_{i,j}$ is
\begin{eqnarray*}
\a(T_{i,j}) &=&  \frac{1}{2N^4}\Bigg|
ij \Bigg( \Big( N-(i+1) \Big) \Big( N-j \Big)-
\Big( N-i \Big) \Big( N- (j+1)\Big) \Bigg) \\
&& + (i+1)j \Bigg( \Big( N-i\Big) \Big( N-\left(j+1 \right) \Big)-
\Big( N-i \Big) \Big( N- j\Big) \Bigg) \\
&& + i(j+1) \Bigg( \Big( N-i\Big) \Big( N-j \Big)-
\Big( N-(i+1) \Big) \Big( N- j\Big) \Bigg)
\Bigg| \\
&=& \frac{1}{2N^3}\big|i-j\big| = \frac{1}{2N^3}\big(j - i\big).
\end{eqnarray*}

The expression for $\a(T'_{i,j})$ can be obtained quickly from the expression for $\a(T_{i,j})$ by swapping $i$ with $i+1$, and $j$ with $j+1$.
Hence
$$\a(T'_{i,j}) = \frac{1}{2N^3}\big|(i+1) - (j+1) \big|=\frac{1}{2N^3}\big(j - i\big).$$
Hence triangles $T_{i,j}$ and $T'_{i,j}$ indeed have equal areas. Could they be congruent? Let $d_{i,j}$ denote the length of the side that triangles $T_{i,j}$ and $T'_{i,j}$ have in common (dotted line in Figure~\ref{fig:quadri1}). This is the line segment with endpoints $P_{i,j+1}$ and $P_{i+1,j}$, and is also a diagonal of the quadrangle $Q_{i,j}$. Then $T_{i,j}$ is a triangle with sides of lengths $s_i$, $s_j$, and $d_{i,j}$, and $T'_{i,j}$ is a triangle with sides of lengths $s_{i+1}$, $s_{j+1}$, and $d_{i,j}$. If triangles $T_{i,j}$ and $T'_{i,j}$ are congruent, then $\{ s_i, s_j \}=\{s_{i+1},s_{j+1} \}$, and hence by our previous observation on these line segments and the assumption $i<i+1<j<j+1$, we must have $\{s_{i+1},s_{j+1} \}=\{ s_{N-i},s_{N-j}\}$. In fact, since $i+1<j+1$ and $N-j<N-i$, it must be that $i+1=N-j$. Thus $T_{i,j}$ and $T'_{i,j}$ can be congruent, but only in the special quadrangle $Q_{i,N-i-1}$; in general, they are incongruent.

However, triangle $T'_{i,j}$ is of course congruent to triangle $T_{i+1,j+1}$ since they both have sides of lengths $s_{i+1}$ and $s_{j+1}$ enclosing the same angle. Consequently, their third sides are equal as well; that is, $d_{i,j}=d_{i+1,j+1}$. That means that quadrilaterals $Q_{i,j}$ and $Q_{i+1,j+1}$ have their lower-left-to-upper-right diagonals of the same length! As we shall see, this is just one of the many properties that quadrilaterals on the same diagonal  have in common.

\begin{figure}[t]
\centering
\includegraphics[scale=0.45]{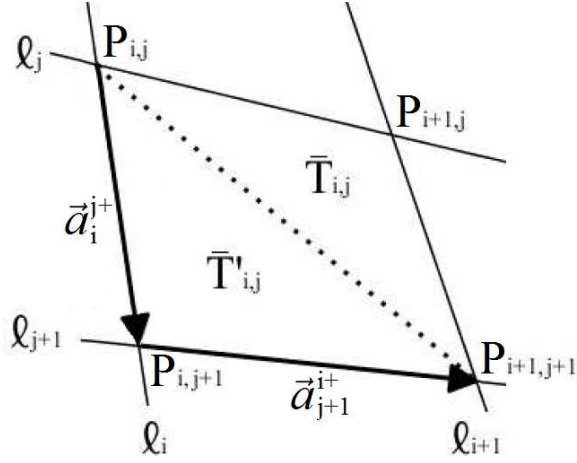}
\caption{Quadrilateral $Q_{i,j}$, and triangles $\bar{T}_{i,j}$ and $\bar{T}'_{i,j}$.}
\label{fig:quadri2}
\end{figure}

Now, let's have a look at what happens when the quadrilateral is cut along its other diagonal, that is, from its top left corner to its bottom right corner (Figure~\ref{fig:quadri2}).
Let $\bar{T}_{i,j}$ denote the triangle with vertices $P_{i,j}, P_{i+1,j}$, and $P_{i+1,j+1}$, and let  $\bar{T'}_{i,j}$ denote the triangle with vertices $P_{i,j}, P_{i,j+1}$, and $P_{i+1,j+1}$. The area of the first triangle is

\begin{eqnarray*}
\a(\bar{T}_{i,j}) &=&  \frac{1}{2N^4}\Bigg|
ij \Bigg( \Big( N-(i+1) \Big) \Big( N-j \Big)-
\Big( N-(i+1) \Big) \Big( N- (j+1)\Big) \Bigg) \\
&& + (i+1)j \Bigg( \Big( N-(i+1)\Big) \Big( N-\left(j+1 \right) \Big)-
\Big( N-i \Big) \Big( N- j\Big) \Bigg) \\
&& + (i+1)(j+1) \Bigg( \Big( N-i\Big) \Big( N-j \Big)-
\Big( N-(i+1) \Big) \Big( N- j\Big) \Bigg)
\Bigg| \\
&=&  \frac{1}{2N^3}\big|i-j+1\big|=\frac{1}{2N^3}\big(j - i - 1\big).
\end{eqnarray*}
The expression for $\a(\bar{T'}_{i,j})$ is obtained from the expression for $\a(\bar{T}_{i,j})$ by interchanging the roles of $i$ and $j$. Hence we immediately obtain
$$\a(\bar{T'}_{i,j}) = \frac{1}{2N^3}\big|j-i+1\big|=\frac{1}{2N^3}(j - i + 1).$$
Observe that $\a(\bar{T}_{i,j}) \ne \a(\bar{T'}_{i,j})$.
Thus, we can see that the way in which the  quadrilateral $Q_{i,j}$ is divided into two triangles determines whether or not the resulting triangles are equal in area.

Using the observation $\a(T_{i,j})=\a(T'_{i,j})$, the area of the quadrilateral $Q_{i,j}$ itself is now easily computed as
$$\a(Q_{i,j})=2 \cdot \a(T_{i,j}) = \frac{1}{N^3}(j - i).$$
Note that this expression depends only on $j-i$, so here is our next surprise:  the quadrilaterals $Q_{i,i+d}$, for all $i \in \{0,1,\ldots,N-3\}$ and all $d\in \{ 2,\ldots, N-i-1\}$ --- these are precisely the quadrilaterals that appear on the same diagonal of the net  --- have equal areas, even though in general they are certainly not congruent! (We invite the readers to verify the latter claim.)

Are there any such surprises in stock for the triangles on the upper-right boundary of our string net? Let's have a look.
Let $T_{i}$ denote the triangle with vertices $P_{i,i+1}$, $P_{i,i+2}$, and $P_{i+1,i+2}$. Again, we use Formula~(\ref{area}) to determine the areas of these triangles.

\begin{eqnarray*}
\a(T_{i}) &=&  \frac{1}{2N^4}\Bigg|
i(i+1) \Bigg( \Big( N-i \Big) \Big( N-(i+2) \Big)-
\Big( N-(i+1) \Big) \Big( N- (i+2)\Big) \Bigg) \\
&& + i(i+2) \Bigg( \Big( N-(i+1)\Big) \Big( N-\left(i+2 \right) \Big)-
\Big( N-i \Big) \Big( N- (i+1)\Big) \Bigg) \\
&& + (i+1)(i+2) \Bigg( \Big( N-i\Big) \Big( N-(i+1) \Big)-
\Big( N-i \Big) \Big( N- (i+2) \Big) \Bigg)
\Bigg| \\
&=&  \frac{1}{N^3}.
\end{eqnarray*}
Observe that the areas of triangles $T_{i}$ do not depend on $i$ at all! Hence, just like quadrilaterals on the same diagonal, all of these triangles, which are also lined up on a ``diagonal'' in Figures~\ref{fig:net} and \ref{fig:checkerednet}, have equal areas although they are in general not congruent.

\bigskip

To summarize, the string art net with a right-angled frame and equidistant intercepts  exhibits the following interesting characteristics.
\begin{description}
\item[{\bf (C1)}] On each line $\ell_i$, for $1 \le i \le N-1$, the intersection points with other lines are equally spaced, with the exception of (consecutive) points $P_{i,i-1}$ and $P_{i,i+1}$, which occur at twice the distance between any other consecutive points. In other words,
$$\vec{a}_{i}^{0+}= \vec{a}_{i}^{1+}=\ldots= \vec{a}_{i}^{(i-2)+}=\frac{1}{2}\vec{b}_{i}= \vec{a}_{i}^{(i+1)+}=\ldots= \vec{a}_{i}^{(N-1)+}.$$
\item[{\bf (C2)}] Each quadrilateral $Q_{i,j}$ is divided into triangles $T_{i,j}$ and $T'_{i,j}$ that are in general incongruent but have equal areas.
\item[{\bf (C3)}] For a fixed $d\in \{ 2,3, \ldots, N-1\}$, the quadrilaterals $Q_{i,i+d}$,  for all $i \in \{ 0,1,\ldots, N-d-1\}$,  have equal areas but are not necessarily congruent.
\item[{\bf (C4)}] Triangles $T_i$, for all $i \in \{0,1,\ldots,N-2\}$, are also have equal areas but are not necessarily congruent.
\end{description}

Having evaluated the areas of the quadrilaterals $Q_{i,j}$ and triangles $T_i$, Properties {\bf (C3)} and {\bf (C4)}  allow us to easily calculate the total area under the piecewise-linear curve.

First, the sum of the areas of the triangles $T_i$ is
$$
\sum_{i=0}^{N-2} \a(T_i) = \sum_{i=0}^{N-2} \frac{1}{N^3} = \frac{1}{N^3} (N-1).$$

Next, we sum up the areas of all quadrilaterals $Q_{i,j}$, taking advantage of the property that $\a(Q_{i,j})$ depends only on $j-i$, and using Formulas (\ref{sum}) and (\ref{squares}).
\begin{eqnarray*}
\sum_{i=0}^{N-3} \sum_{j=i+2}^{N-1} \a(Q_{i,j}) &=& \sum_{i=0}^{N-3} \sum_{d=2}^{N-i-1} \a(Q_{i,i+d}) \\
&=& \sum_{i=0}^{N-3} \sum_{d=2}^{N-i-1} \frac{1}{N^3}d =
\frac{1}{N^3} \sum_{i=0}^{N-3} \sum_{d=2}^{N-i-1} d \\
&=& \frac{1}{N^3} \sum_{i=0}^{N-3} \Bigg( \frac{1}{2}(N-i-1)(N-i) -1\Bigg)\\
&=& \frac{1}{2N^3} \sum_{i=0}^{N-3} \Bigg( (N-i)^2-(N-i)-2 \Bigg)\\
&=& \frac{1}{2N^3} \sum_{i=3}^{N} \Bigg( i^2-i-2 \Bigg) \\
 &=&\frac{1}{2N^3} \Bigg( \Big( \frac{1}{6}N(N+1)(2N+1)-(1+2^2) \Big)  - \Big( \frac{1}{2}N(N+1)-(1+2) \Big) - 2(N-2) \Bigg) \\
&=& \frac{1}{6N^3}(N-1)(N-2)(N+3)
\end{eqnarray*}
Hence the area of the region in the first quadrant under the piecewise-linear curve determined by the points $P_{0,1},P_{1,2},\ldots,P_{N-1,N}$ is
$$A=\frac{1}{N^3} (N-1)+\frac{1}{6N^3}(N-1)(N-2)(N+3)= \frac{1}{6N^2}(N - 1)(N + 1).$$
Observe that as $N \rightarrow \infty$, the area of this region approaches $\frac{1}{6}$. This is --- as it should be --- precisely the area in the first quadrant under the envelope of our family of lines. Namely, it was shown in \cite{Que} that this envelope has the equation $y=x-2x^{\frac{1}{2}}+1$, and with a little bit of calculus we obtain the area under the curve as
$$\int_0^1 \Big(x-2x^{\frac{1}{2}}+1 \Big)dx=\Big(\frac{1}{2}x^2-\frac{4}{3}x^{\frac{3}{2}}+x \Big) \Big\vert_0^1=\frac{1}{6}.$$

\bigskip

Having discovered surprising properties {\bf (C1)}--{\bf (C4)} of the string art net with a right-angled frame and equidistant pegs, it is natural to ask whether these properties are preserved when the frame is not right-angled, or when the pegs are not equidistant. We shall explore these two questions in the next two sections.

\section{Skew frame}

\begin{figure}[t]
\centering
\includegraphics[scale=0.8]{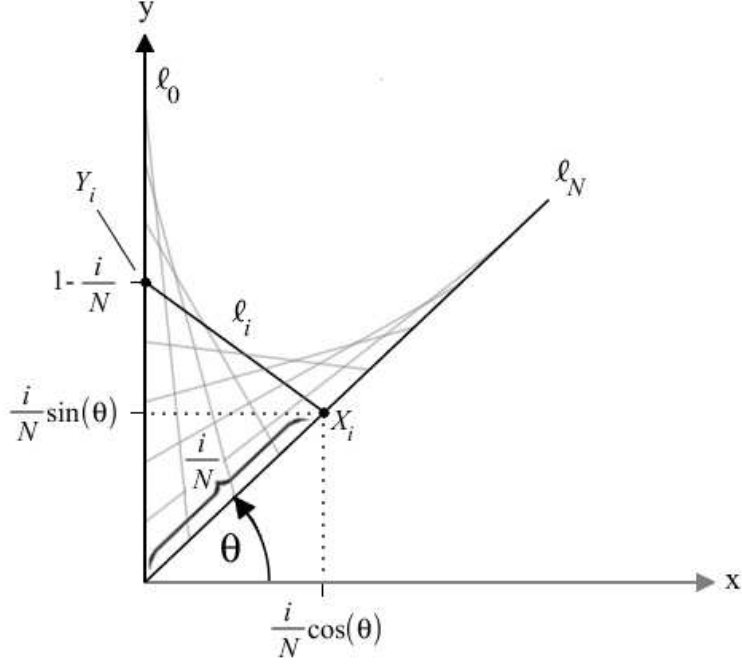}
\caption{A string art net with an acute-angled frame.}
\label{fig:skewed}
\end{figure}

In this section, we again assume that the pegs on both axes are equidistant, but we shall slant the previously horizontal axis of the frame at an angle  $\theta$, for $-\frac{\pi}{2} < \theta < \frac{\pi}{2}$, with respect to the $x$-axis of our coordinate system. Note that when $0< \theta < \frac{\pi}{2}$, the resulting net has an acute-angled frame (shown in Figure~\ref{fig:skewed}), while for $-\frac{\pi}{2} < \theta < 0$, the frame is obtuse-angled. Our calculations, however, will apply to both cases.

Denote $a = \cos (\theta) $ and  $b = \sin (\theta)$. As in the previous case, our string art net will be framed by lines $\ell_0$ and $\ell_N$, however, line $\ell_N$ will be rotated by the angle $\theta$. In other words, $\ell_N$ is the line with equation $y=\frac{b}{a}x$, and for all $i \in I$, line  $\ell_i$ passes through the points
$$X_i = \left(\frac{i}{N}a,\frac{i}{N}b \right) \qquad \mbox{and} \qquad Y_i = \left(0,1-\frac{i}{N}\right).$$

Similarly to our calculation in Section 3, the equation of line $\ell_i$ is obtained as
$$y = \frac{(1-\frac{i}{N}) - \frac{i}{N}b}{0-\frac{i}{N}a}\left(x - \frac{i}{N}a\right) + \frac{i}{N}b$$
with canonical form
$$y = -\frac{N - i(b + 1)}{ai}x + \Big(1 - \frac{i}{N}\Big). $$

The $x$-coordinate of the intersection point of lines $\ell_i$ and $\ell_j$ is the solution to the equation
$$-\frac{N - i(b + 1)}{ai}x + \Big(1 - \frac{i}{N}\Big) = -\frac{N - j(b + 1)}{aj}x + \Big(1 - \frac{j}{N}\Big),$$
that is,
$$ x = \frac{a}{N^2}ij$$
The $y$-coordinate of the intersection point is
$$y = -\frac{N - i(b + 1)}{ai} \cdot \frac{a}{N^2}ij + \Big(1 - \frac{i}{N}\Big),$$
that is,
$$ y = \frac{1}{N^2}\Big( (N-i)(N-j) + bij \Big).$$
Therefore, the intersection of lines  $\ell_i$ and $\ell_j$ (for $i \ne j$) is the point
$$P_{i,j}=\frac{1}{N^2}\Big(aij,(N-i)(N-j) + bij \Big).$$

To examine Properties {\bf (C1)}--{\bf (C4)}, as in Section 3, we next define  vectors $\vec{a}_{i}^{j+}$ and $\vec{b}_{i}$. For all $i, j \in I$ such that   $j \not\in \{i-1,i,N\}$, let
$$\vec{a}_{i}^{j+} = P_{i,j+1} - P_{i,j} =
\frac{1}{N^2}\Bigg( ai, -(N-i)+bi \Bigg),$$
and for all $i \in \{1,2,\ldots, N-1\}$, let
$$\vec{b}_{i} = P_{i,i+1} - P_{i,i-1} = \frac{2}{N^2}\Bigg( ai, -(N-i)+bi \Big) \Bigg).$$
Interestingly, the vectors exhibit Property {\bf (C1)} as in the case of the right-angled frame: vector $\vec{a}_{i}^{j+}$ is independent of the parameter $j$, and scaled by 2 results in vector $\vec{b}_{i}$. This observation suggests that the string art net will display the remaining interesting properties observed in Section 3, despite the different angle of the frame. Let's verify this conjecture.

Again, denote the quadrilateral with vertices $P_{i,j}, P_{i+1,j}$,  $P_{i,j+1}$, and $P_{i+1,j+1}$ by $Q_{i,j}$, and partition each quadrilateral in two different ways, as in Figures~\ref{fig:quadri1} and \ref{fig:quadri2}.
As before, let  $T_{i,j}$ denote the triangle with vertices $P_{i,j}, P_{i+1,j}$, and $P_{i,j+1}$, and let  $T'_{i,j}$ denote the triangle with vertices $P_{i+1,j}, P_{i,j+1}$, and $P_{i+1,j+1}$ (see Figure~\ref{fig:quadri1}). The area of  triangle $T_{i,j}$ is then
\begin{eqnarray*}
\a(T_{i,j}) &=&  \frac{a}{2N^4}\Bigg|
ij \Bigg( \Big(N-(i+1) \Big) \Big( N-j \Big)+b(i+1)j -
\Big(N-i \Big) \Big( N-(j+1) \Big) - bi(j+1) \Bigg) \\
&& + (i+1)j \Bigg(
\Big(N-i \Big) \Big( N-(j+1) \Big) + bi(j+1)
-\Big(N-i \Big) \Big( N-j \Big)-bij \Bigg) \\
&& + i(j+1) \Bigg(  \Big(N-i \Big) \Big( N-j \Big)+bij
-\Big(N-(i+1) \Big) \Big( N-j \Big)-b(i+1)j \Bigg)
\Bigg| \\
&=&  \frac{a}{2N^3}\big(j - i\big).
\end{eqnarray*}
(Note that because of the analogous calculation performed in Section 3, we only needed to show here that all terms with a factor $b$ cancel each other out.)
Similarly, we obtain
$$\a(T'_{i,j}) = \frac{a}{2N^3}(j - i).$$
We see once again that quadrilateral $Q_{i,j}$ is divided into two in general incongruent triangles $T_{i,j}$ and $T'_{i,j}$ of equal areas, so Property {\bf (C2)} holds.

On the other hand, let $\bar{T}_{i,j}$ denote the triangle with vertices $P_{i,j}, P_{i+1,j}$, and $P_{i+1,j+1}$, and let  $\bar{T'}_{i,j}$ denote the triangle with vertices $P_{i,j}, P_{i,j+1}$, and $P_{i+1,j+1}$ (Figure~\ref{fig:quadri2}). The areas of these triangles are then

\begin{eqnarray*}
\a(\bar{T}_{i,j}) &=&  \frac{a}{2N^4}\Bigg|
ij \Bigg(
 \Big(N-(i+1) \Big) \Big( N-j \Big) + b(i+1)j
 -\Big(N-(i+1) \Big) \Big( N-(j+1) \Big) - b(i+1)(j+1)
\Bigg) \\
&& + (i+1)j \Bigg(
 \Big(N-(i+1) \Big) \Big( N-(j+1) \Big) + b(i+1)(j+1)
 -\Big(N-i \Big) \Big( N-j \Big) - bij
\Bigg) \\
&& + (i+1)(j+1) \Bigg(
\Big(N-i \Big) \Big( N-j \Big) + bij
-\Big(N-(i+1) \Big) \Big( N-j \Big) - b(i+1)j
 \Bigg)
\Bigg| \\
&=&  \frac{a}{2N^3}\big( j-i-1\big)
\end{eqnarray*}
and
$$\a(\bar{T'}_{i,j}) = \frac{a}{2N^3}\big( j-i+1\big),$$
so we can see that $\a(\bar{T}_{i,j}) \ne \a(\bar{T'}_{i,j})$.
As expected, the way in which the quadrilateral is divided determines  whether the two parts have equal areas or not, regardless of the frame angle of the string art net.

Next, using Property {\bf (C2)}, we find the area of the quadrilateral $Q_{i,j}$ as
$$\a(Q_{i,j}) = 2 \cdot \a(T_{i,j})=\frac{a}{N^3}(j - i)$$
--- an expression dependent only on $j-i$,
and conclude that Property {\bf (C3)} holds as well.

Finally, let $T_{i}$ be the triangle with vertices $P_{i,i+1}$, $P_{i,i+2}$, and $P_{i+1,i+2}$. Its area is obtained as

\begin{eqnarray*}
\a(T_{i}) &=&  \frac{a}{2N^4}\Bigg|
i(i+1) \Bigg(
\Big( N-i \Big) \Big( N-(i+2) \Big)+ bi(i+2)
- \Big( N-(i+1) \Big) \Big( N-(i+2) \Big)- b(i+1)(i+2)
\Bigg) \\
&& + i(i+2)\Bigg(
\Big( N-(i+1) \Big) \Big( N-(i+2) \Big)+ b(i+1)(i+2)
-\Big( N-i \Big) \Big( N-(i+1) \Big)- bi(i+1)
 \Bigg) \\
&& + (i+1)(i+2) \Bigg(
\Big( N-i \Big) \Big( N-(i+1) \Big)+ bi(i+1)
-\Big( N-i \Big) \Big( N-(i+2) \Big)- bi(i+2)
 \Bigg)
\Bigg| \\
&=&  \frac{a}{N^3},
\end{eqnarray*}
which is again independent of $i$, so we may conclude that Property {\bf (C4)} is satisfied as well. Thus, as predicted, the string art net in the skewed case displays the same interesting properties {\bf (C1)}--{\bf (C4)} as in the right-angled case, regardless of the angle $\theta$.

As a bonus, observe that the areas of triangles $T_{i,j}$, $T'_{i,j}$, and $T_i$ in this case are equal to the areas of the corresponding triangles in Section 3 multiplied by $a$. Therefore, we can immediately conclude that the area of the region bounded by the piecewise-linear curve through the points $P_{0,1},P_{1,2},\ldots,P_{N-1,N}$ and lines $\ell_0$ and $\ell_N$ is
$$A = \frac{a}{6N^2}(N - 1)(N + 1).$$
Observe that, surprisingly, since $a=\cos \theta = \cos (- \theta)$,  the area of the region of  a string art net with an acute angle and of a string art net with an obtuse angle --- which have very different shapes! --- will be equal as long as the angles between  line $\ell_N$ and the $x$-axis are equal in absolute value.

\section{Changing the spacing function}

In this section, we are asking the question whether, assuming the right-angled frame, any other {\em spacing function} (as named in \cite{Que}), which determines the placement of the pegs along the frame,  yields the same fascinating Properties {\bf (C1)}--{\bf (C4)} of the string art net that we observed in Sections 3 and 4.

Our spacing function, which we call $f$, will give the $x$-intercept of line $\ell_i$. In Section 3 we had $f(i)=\frac{i}{N}$ so that line $\ell_i$ had the intercepts\footnote{Observe that (\ref{intercepts-basic}) can also be generalized as $X_i = \Big( f(i),0 \Big)$ and $Y_i = \Big( 0,f(N-i) \Big)$, but this setting does not lend itself to an accessible analysis.}
\begin{eqnarray}
X_i = \Big( f(i),0 \Big) & \mbox{and} & Y_i = \Big( 0,1-f(i) \Big).
\label{intercepts}
\end{eqnarray}
We shall now vary function $f$ but assume the same intercepts of our lines $\ell_i$ relative to $f$. Note that our spacing function should satisfy the following conditions:
\begin{enumerate}[(i)]
\item $f$ is strictly increasing,
\item $f: I \rightarrow [0,1]$, and
\item $f(0) = 0$ and $f(N)= 1$,
\end{enumerate}
so that our string art net remains enclosed in the unit square $[0,1] \times [0,1]$.
For better readability, we shall henceforth denote $f(i)=x_i$, for all $i \in I$. Our analysis below will be most relevant for $N \ge 4$, but we shall comment on the cases with $N \le 3$ at the end.

Observe that line $\ell_i$ is then described by the equation
$$y = \frac{x_{i}-1}{x_{i}}x + 1-x_{i}.$$
The intersection of lines $\ell_i$ and $\ell_j$, for $i,j \in I$ such that $i \ne j$, is the point
$$P_{i,j}=\Big( x_{i}x_{j},(x_{i}-1)(x_{j}-1) \Big),$$
and the vectors between two consecutive points on line $\ell_i$ are
$$\vec{a}_{i}^{j+} = P_{i,j+1} - P_{i,j} = ( x_{j+1}-x_{j}) \Big( x_{i}, x_{i}-1 \Big)$$
for $j \not\in \{i-1,i,N\}$, and
$$\vec{b}_{i} = P_{i,i+1} - P_{i,i-1} = ( x_{i+1}-x_{i-1}) \Big( x_{i},x_{i}-1 \Big)$$
for $i \in \{1,\ldots,N-1\}$.
Observe that in general, the length of vector $\vec{a}_{i}^{j+}$ depends on parameter $j$. We would like to determine all functions $f$ --- or equivalently, sequences $(x_{0},\ldots,x_{N})$ --- for which vector $\vec{a}_{i}^{j+}$ depends only on the parameter $i$;  this is the nice Property {\bf (C1)} we observed in Section 3. It is clear that, in order to obtain this property,  $x_{j+1}-x_{j}$ must be constant. Since $x_{0}=0$ and $x_{N}=1$, it then follows that $x_{j}=\frac{j}{N}$, so $f$ is precisely the spacing function used in Section 3.

It appears that we have encountered a dead end right from the start. But is it possible that some of the other interesting properties observed in Section 3 can be satisfied without $\vec{a}_{i}^{j+}$ being independent from $j$? In other words, could some of Properties {\bf (C2)}--{\bf (C4)} hold even though {\bf (C1)} does not?

Let symbols $Q_{i,j}$, $T_{i,j}$, $T'_{i,j}$, and $T_i$  denote the quadrilaterals and triangles  as defined in Section 3.
To examine Property {\bf (C2)}, determine the areas of triangles $T_{i,j}$ and $T'_{i,j}$. Assuming $0\le i<i+1<j<j+1 \le N$, we obtain

\begin{eqnarray*}
\a(T_{i,j}) &=&  \frac{1}{2}\Bigg|
x_{i}x_{j} \Big( (x_{i+1}-1) (x_{j}-1)- (x_{i}-1) (x_{j+1}-1)\Big) \\
&& + x_{i+1}x_{j} \Big( ( x_{i}-1) (x_{j+1}-1)- (x_{i}-1) (x_{j}-1)\Big) \\
&& + x_{i}x_{j+1} \Big( (x_{i}-1) ( x_{j}-1) - (x_{i+1}-1) (x_{j}-1) \Big)
\Bigg| \\
&=&  \frac{1}{2} (x_{i+1}-x_{i})(x_{j}-x_{i}) (x_{j+1}-x_{j}),
\end{eqnarray*}
and, swapping $i$ with $i+1$, and $j$ with $j+1$,
$$\a(T'_{i,j}) = \frac{1}{2}  (x_{i+1}-x_{i})(x_{j+1}-x_{i+1}) (x_{j+1}-x_{j}).$$
Observe that $\a(T_{i,j})=\a(T'_{i,j})$ implies
$$x_{j}-x_{i}=x_{j+1}-x_{i+1},$$
in other words,
$$x_{i+1}-x_{i}=x_{j+1}-x_{j}.$$
If this is true for all $i$ and $j$, then  $x_{i+1}-x_{i}$
is constant. We deduce that for Property {\bf (C2)} to hold, our spacing function must again be defined as $x_{i}=\frac{i}{N}$.

Let us now examine Properties {\bf (C3)} and {\bf (C4)}.
For any spacing function $f$, we obtain that the area of the quadrilateral $Q_{i,j}$ is
\begin{eqnarray}
\a(Q_{i,j})=\a(T_{i,j})+\a(T'_{i,j})=
\frac{1}{2}(x_{i+1}-x_{i})(x_{j+1}-x_{j})(x_{j}+x_{j+1}-x_{i}-x_{i+1}), \label{Qij}
\end{eqnarray}
and the area of triangle $T_i$ is
\begin{eqnarray}
\a(T_{i}) &=&
\frac{1}{2}\Bigg|
x_{i}x_{i+1} \Big( (x_{i}-1) (x_{i+2}-1)- (x_{i+1}-1) (x_{i+2}-1)\Big) \nonumber \\
&& + x_{i}x_{i+2} \Big( ( x_{i+1}-1) (x_{i+2}-1)- (x_{i}-1) (x_{i+1}-1)\Big) \nonumber \\
&& + x_{i+1}x_{i+2} \Big( (x_{i}-1) ( x_{i+1}-1) - (x_{i}-1) (x_{i+2}-1) \Big)
\Bigg| \nonumber \\
&=& \frac{1}{2} (x_{i+1}-x_{i})(x_{i+2}-x_{i})(x_{i+2}-x_{i+1}). \label{Ti}
\end{eqnarray}

If quadrilaterals and triangles on the same diagonal of the net are to have equal areas, then $\a(T_{i})$ should be constant, and  $\a(Q_{i,j})$ should depend only on $j-i$. Is this possible without having the vectors $\vec{a}_{i}^{j+}$ depend only on $i$, in other words, without $f$ being defined by $x_{i}=\frac{i}{N}$?

Assuming  $N \ge 4$, we have $N-2$ equations (for $i=0,1,\ldots,N-3$) of the form
$$\a(T_i)=\a(T_{i+1}),$$
together with equation $$\a(Q_{0,2})=\a(Q_{1,3}),$$
which gives $N-1$ equations for the variables $x_{1}, x_{2}, \ldots, x_{N-1}$. As we show below, this system --- with our constraints (i)-(iii) --- indeed has a unique solution.

To formulate a hypothesis about the $x$-intercepts $x_{i}$,  we first obtain the expressions for $x_{3}$, $x_{4}$, and $x_5$. From
$$\a(T_0) = \a(T_1), \label{T0-1}$$
using (\ref{Ti}), we obtain
$$(x_{1}-x_{0})(x_{2}-x_{0})(x_{2}-x_{1}) = (x_{2}-x_{1})(x_{3}-x_{1})(x_{3}-x_{2}),$$
which, since $x_0=0$, simplifies to
$$x_{3} \Big( x_{3} - (x_{1} + x_{2}) \Big) = 0.$$
Hence
$$x_{3} = x_{1}+x_{2}.$$
From here, using expression (\ref{Ti}), equalities $\a(T_1)=\a(T_{2})$ and $\a(T_2)=\a(T_{3})$ yield equations
$$x_{4} = 2x_{2} \quad \mbox{ and } \quad x_{5} = x_{1}+2x_{2},$$
respectively.
We thus hypothesize that the general expression for $x_{i}$ is as follows.

\smallskip

\begin{claim} For $i=1,\ldots,N$,
\begin{eqnarray}
x_{i}=\left\{ \begin{array}{ll}
              x_{1}+\frac{i-1}{2}x_{2} & \mbox{ if } i \mbox{ is odd} \\
              \frac{i}{2}x_{2} & \mbox{ if } i \mbox{ is even} \\
              \end{array}\right. .
              \label{piecewise}
\end{eqnarray}
\end{claim}

\bigskip

\noindent
To prove this claim, we shall use strong induction.

We have seen already that  (\ref{piecewise}) is true for $1 \le i \le 5$, so in particular, the basis of induction holds.

Now assume that for some $k \in \{2,\ldots,N-1\}$, Formula~(\ref{piecewise}) holds for all $i \leq k$. This is our induction hypothesis. We need to show that then (\ref{piecewise}) also holds for $i = k + 1$.

Observe that, since $k \geq 2$,
$$\a(T_{k-2})=\a(T_{k-1}),$$
which, using (\ref{Ti}), simplifies to
$$(x_{k-1}-x_{k-2})(x_{k}-x_{k-2})=(x_{k+1}-x_{k-1})(x_{k+1}-x_{k}).$$
Rearranging the terms, we obtain
$$(x_{k-1}+x_{k})(x_{k+1}-x_{k-2}) = (x_{k-2}+x_{k+1})(x_{k+1}-x_{k-2}).$$
Since $f$ is strictly increasing, we conclude that
$$x_{k-1}+x_{k} =x_{k-2}+x_{k+1},$$
that is,
$$x_{k+1} = x_{k-1}+x_{k}-x_{k-2}.$$
We shall use this recurrence relation to show that (\ref{piecewise}) for $i=k+1$ follows from our induction hypothesis.

If $k$ is odd, the induction hypothesis yields
$$x_{k} = x_{1} + \frac{k-1}{2}x_{2}, \quad x_{k-1} = \frac{k-1}{2}x_{2}, \quad \mbox{ and } \quad x_{k-2} = x_{1} + \frac{k-3}{2}x_{2},$$
whence
$$x_{k+1} = x_{k-1}+x_{k}-x_{k-2} = \frac{k+1}{2}x_{2}$$
as required.

Similarly, if $k$ is even,  by the induction hypothesis
$$x_{k} = \frac{k}{2}x_{2}, \quad
x_{k-1} = x_{1} + \frac{k-2}{2}x_{2}, \quad \mbox{ and } \quad
x_{k-2} = \frac{k-2}{2}x_{2},$$
yielding
$$x_{k+1} = x_{k-1}+x_{k}-x_{k-2} = x_{1} + \frac{k}{2}x_{2}.$$
Hence Formula~(\ref{piecewise})  for $i=k+1$ follows from our induction hypothesis, and the induction step is complete. By the Principle of Mathematical Induction,  our Claim follows.

It remains to find the values of $x_1$ and $x_2$. We  use
$$\a(Q_{0,2})=\a(Q_{1,3})$$
together with (\ref{Qij}) to obtain
\begin{eqnarray}
\frac{1}{2}(x_1-x_0) (x_{3}-x_{2})(x_{2}+x_{3}-x_{0}-x_1) &=& \frac{1}{2}(x_{2}-x_{1})(x_{4}-x_{3})(x_{3}+x_{4}-x_{1}-x_{2}),
\label{lasteq}
\end{eqnarray}
which will give the first equation relating $x_1$ and $x_2$. Namely, substituting
$$x_0=0, \quad x_3=x_1 + x_2, \quad \mbox{ and } \quad x_4=2x_2,$$
(\ref{lasteq}) reduces to
$$x_{2} = 2x_{1}.$$

The second equation relating $x_1$ and $x_2$ is
$$1=
\left\{ \begin{array}{ll}
                 x_{1}+\frac{N-1}{2}x_{2} & \mbox{ if } N \mbox{ is odd} \\
                  \frac{N}{2}x_{2}& \mbox{ if } N \mbox{ is even}
               \end{array}\right.,
$$
obtained from the assumption $x_N=1$ and (\ref{piecewise}). Combining this result with $x_{2} = 2x_{1}$ then yields
$$x_{1}=\frac{1}{N} \quad \mbox{ and } \quad x_{2}=\frac{2}{N}.$$
Finally, substituting these expressions in (\ref{piecewise}) we obtain $x_i=\frac{i}{N}$.

We conclude that, for $N \ge 4$,  the only spacing function satisfying assumptions (i)-(iii) that gives rise to a string art net with intercepts (\ref{intercepts}) and Properties {\bf (C3)} and {\bf (C4)} is indeed the equidistant spacing function $f(i)=\frac{i}{N}$.

How about small nets with $N \le 3$? For $N=1$, our net is reduced to lines $\ell_0$ and $\ell_N$, and for $N=2$, it contains a single triangle and no quadrilaterals, so Properties {\bf (C1)}--{\bf (C4)} hold vacuously. For $N=3$, the net consists of two triangles and one quadrilateral. From (\ref{T0-1}) we obtain
$$1=x_3=x_1+x_2,$$
so $x_2=1-x_1$. Thus line $\ell_1$ has intercepts $(x_1,0)$ and $(0,1-x_1)$, line $\ell_2$ has intercepts $(1-x_1,0)$ and $(0,x_1)$, and our net is symmetric about the axis $y=x$. For $f$ to be strictly increasing, we need $0<x_1 <\frac{1}{2}$, however, Properties {\bf (C2)}--{\bf (C4)} hold for any $x_1$ in this range! To the contrary, Property {\bf (C1)} requires that $\vec{b}_1=2\vec{a}_1^{2+}$, which can be seen to imply $x_1=\frac{1}{3}$, again resulting in the equidistant spacing function $f(i)=\frac{i}{N}$.

\section{Summary}

We showed that the string art net of the type illustrated in Figure~\ref{fig:art} exhibits unexpected symmetry.
In particular, we showed that each line of the net is divided by other lines into segments, all  but one of the same length and one of twice the length of the others. Furthermore, quadrilaterals in the net that are arranged along a diagonal from the upper left to the lower right corner are pairwise incongruent but have equal areas, and the same is true of triangles formed along the upper-right border of the net. Finally, we showed that these properties are preserved when the angle between the axes changes, but hold only when consecutive pegs on the axes are equidistant.

\end{document}